\documentclass[english, 11pt]{amsart}

\usepackage{babel}
\usepackage[utf8x]{inputenc}
\usepackage{mathrsfs}
\usepackage{amsmath}
\usepackage{amssymb}
\usepackage{amsthm}
\usepackage[all, cmtip]{xy}
\usepackage{leftidx}
\usepackage{amsmath,amssymb,amsfonts}
\usepackage{amsmath,amscd}
\usepackage{textcomp}
\usepackage{pictexwd,dcpic}

\newtheorem{thm}{Theorem}[section]
\newtheorem{lemma}[thm]{Lemma}
\newtheorem{remark}[thm]{Remark}
\newtheorem{proposition}[thm]{Proposition}
\newtheorem{question}[thm]{Question}

\newtheorem{definition}[thm]{Definition}

\theoremstyle{definition}

\newcommand{\dbar}{\overline{\partial}}
\newcommand{\ddbar}{\partial \overline{\partial}}

\newcommand\ac{{\rm ac}}
\newcommand\reg{{\rm reg}}

\numberwithin{equation}{section}

\begin{document}

\renewcommand{\subjclassname}{%
\textup{2010} Mathematics Subject Classification}

\title[Extension of holomorphic functions and cohomology classes]{Extension of holomorphic functions and cohomology classes from non reduced analytic subvarieties\vskip2mm\strut}

\date{\today}



\thanks{Supported by the European Research Council project
``Algebraic and K\"ahler Geometry'' (ERC-ALKAGE, grant No. 670846
from September 2015).}

\author[Jean-Pierre Demailly]{Jean-Pierre Demailly${}^*$}
\address{
(*) Institut Fourier, Universit\'e Grenoble Alpes, 100 rue des Maths,
38610 Gi\`eres, France\vskip0pt {\it e-mail}\/:
{\tt jean-pierre.demailly@univ-grenoble-alpes.fr}}
\maketitle

{\hfill\emph{\small in honor of Professor Kang-Tae Kim on the occasion of his sixtieth birthday}\hfill}
\vskip3mm\strut  

\begin{abstract}
The goal of this survey is to describe some recent results
concerning the $L^2$ extension of holomorphic sections or cohomology classes
with values in vector bundles satisfying weak semi-positivity
properties. The results presented here are generalized versions of
the Ohsawa-Takegoshi extension theorem, and borrow many
techniques from the long series of papers by T.~Ohsawa. The recent achievement
that we want to point out is that the surjectivity property holds true
for restriction
morphisms to non necessarily reduced subvarieties, provided these are defined
as zero varieties of multiplier ideal sheaves. The new idea involved
to approach the existence problem is to make use of $L^2$ approximation
in the Bochner-Kodaira technique. The extension results hold under
curvature conditions that look pretty optimal. However, a major unsolved
problem is to obtain natural (and hopefully best possible) $L^2$ estimates
for the extension in the case of non reduced subvarieties -- the case
when $Y$ has singularities or several irreducible components is also
a substantial issue.
\end{abstract}

\section{Introduction and main results}

The problem considered in these notes is whether a holomorphic object
$f$ defined on a subvariety $Y$ of a complex manifold $X$ can be
extended as a holomorphic object $F$ of the same nature on the
whole of~$X$. Here, $Y$ is a subvariety
defined as the the zero zet of a non necessarily reduced ideal $\mathcal{I}$
of $\mathcal{O}_X$, the object to extend can be either a section
$f\in H^0(Y,E_{|Y})$ or a cohomology class
$f\in H^q(Y,E_{\|Y})$, and we look for an extension
$F\in H^q(X,E)$, assuming suitable convexity properties of $X$
and~$Y$, suitable $L^2$ conditions for $f$ on~$Y$, and appropriate curvature
positivity hypotheses for the bundle~$E$. When $Y$ is not connected,
this can be also seen as an interpolation problem  -- the situation where
$Y$ is a discrete set is already very interesting.

The prototype of such results is the celebrated $L^2$ extension theorem
of Ohsawa-Takegoshi \cite{OT87}, which deals with the important case when
$X=\Omega\subset\mathbb{C}^n$ is a pseudoconvex open set, and
$Y=\Omega\cap L$ is the intersection of $\Omega$ with a complex affine
linear subspace $L\subset\mathbb{C}^n$. The accompanying $L^2$ estimates
play a very important role in applications, possibly even more than
the qualitative extension theorems by themselves (cf.\ section 4 below).
The related techniques have then been
the subject of many works since 1987, proposing either greater generality
(\cite{Ohs88,Ohs94,Ohs95,Ohs01,Ohs03,Ohs05}, \cite{Man93}, \cite{Dem97},
\cite{Pop05}), alternative proofs (\cite{Ber96}, \cite{Che11}),
improved estimates (\cite{MV07}, \cite{Var10}) or optimal ones
(\cite{Blo13}, \cite{GZ15a}, \cite{BL16}).

In this survey, we mostly follow the lines of our previous papers
\cite{Dem15b} and \cite{CDM17}, whose goal is to pick the
weakest possible curvature and convexity hypotheses, while allowing
the subvariety $Y$ to be non reduced. The ambient complex manifold $X$
is assumed to be a K\"ahler and \emph{holomorphically convex} (and
thus not necessarily compact); by the Remmert reduction theorem, the
holomorphic convexity is equivalent to the existence of a proper
holomorphic map $\pi:X\to S$ onto a Stein complex space $S$, hence
arbitrary relative situations over Stein bases are allowed.
We consider a holomorphic line bundle $E\to X$ equipped with a
singular hermitian metric $h$, namely a metric which can be expressed
locally as $h=e^{-\varphi}$ where $\varphi$ is a \emph{quasi-psh}
function, i.e.\ a function that is locally the sum
$\varphi=\varphi_0+u$ of a plurisubharmonic function $\varphi_0$ and
of a smooth function $u$. Such a bundle admits a curvature current
\begin{equation}
\Theta_{E, h}:=i\ddbar\varphi=i\ddbar\varphi_0+i\ddbar u
\end{equation}
which is locally the sum of a positive $(1,1)$-current 
$i\ddbar \varphi_0$ and a smooth
$(1,1)$-form $i\ddbar u$. Our goal is to extend sections that
are defined on a non necessarily reduced complex subspace $Y\subset X$,
when the structure sheaf 
\hbox{$\mathcal{O}_Y:=\mathcal{O}_X/\mathcal{I}(e^{-\psi})$} is given by the 
multiplier ideal sheaf of a quasi-psh function $\psi$ 
with \emph{neat analytic singularities}, i.e.\
locally on a neighborhood $V$ of an arbitrary point $x_0\in X$ we have
\begin{equation}
\psi(z)=c\log\sum|g_j(z)|^2+v(z),\qquad g_j\in\mathcal{O}_X(V),~~
v\in C^\infty(V).
\end{equation}
Let us recall that the multiplier ideal sheaf $\mathcal{I}(e^{-\varphi})$ of
a quasi-psh function $\varphi$ is defined by
\begin{equation}
\mathcal{I}(e^{-\varphi})_{x_0}=\big\{f\in\mathcal{O}_{X,x_0}\,;\;\exists U\ni x_0\,,\;
\int_U|f|^2e^{-\varphi}d\lambda<+\infty\big\}
\end{equation}
with respect to the Lebesgue measure $\lambda$ in some local coordinates 
near~$x_0$. As is well known,
$\mathcal{I}(e^{-\varphi})\subset\mathcal{O}_X$ is a coherent ideal sheaf
(see e.g.\ \cite{Dem-book}). We also denote by $K_X=\Lambda^nT^*_X$ the 
canonical bundle of an $n$-dimensional complex manifold $X\,$; in the case
of (semi)positive curvature, the Bochner-Kodaira identity yields positive
curvature terms only for $(n,q)$-forms, so the best way to state results
is to consider the adjoint bundle $K_X\otimes E$ rather than the bundle
$E$~itself. The main qualitative statement is given by the following
result of \cite{CDM17}.

\begin{thm}\label{main-theorem}
Let $E$ be a holomorphic line bundle over a holomorphically
convex K\"ahler mani\-fold~$X$. Let $h$ be a possibly singular hermitian 
metric on $E$, $\psi$~a quasi-psh function with neat analytic singularities
on $X$. Assume that 
there exists a positive continuous function $\delta>0$ on $X$ such that
\begin{equation}\label{weak-curv-cond}
\Theta_{E,h}+(1+\alpha\delta)i\ddbar \psi \geq 0\qquad\text{in the sense
of currents, for}~~\alpha=0,1.
\end{equation}
Then the morphism induced by the natural inclusion 
$\mathcal{I}(he^{-\psi}) \to \mathcal{I}(h)$ 
\begin{equation}\label{main-inj}
H^{q}(X, K_X \otimes E \otimes \mathcal{I}(he^{-\psi})) 
\to H^{q}(X, K_X \otimes E \otimes \mathcal{I}(h)) 
\end{equation}
is injective for every $q\geq 0$. 
In other words, the morphism induced by the natural sheaf surjection $\mathcal{I}(h) \to  \mathcal{I}(h)/\mathcal{I}(he^{-\psi})$
\begin{equation}\label{main-surj}
H^{q}(X, K_X \otimes E \otimes \mathcal{I}(h)) \to 
H^{q}(X, K_X \otimes E \otimes \mathcal{I}(h)/\mathcal{I}(he^{-\psi})) 
\end{equation}
is surjective for every $q\geq 0$. 
\end{thm}

\begin{remark}\label{main-consequence} {\rm
{\bf (A)} When $h$ is smooth, we have $\mathcal{I}(h)=\mathcal{O}_X$ and 
$$
\mathcal{I}(h)/\mathcal{I}(he^{-\psi})=\mathcal{O}_X/\mathcal{I}(e^{-\psi})
:=\mathcal{O}_Y
$$
where $Y$ is the zero subvariety of the ideal 
sheaf~$\mathcal{I}(e^{-\psi})$. Hence, the surjectivity statement can
be interpreted an extension theorem 
with respect to the restriction morphism
\begin{equation}\label{main-degree-q}
H^q(X, K_X \otimes E)\to 
H^q(Y, (K_X \otimes E)_{|Y}).
\end{equation}
In general, the quotient sheaf 
$\mathcal{I}(h)/\mathcal{I}(he^{-\psi})$ is supported in
an analytic subvariety $Y\subset X$, which is the zero set of the
conductor ideal
\begin{equation}
\mathcal{J}_Y:=\mathcal{I}(he^{-\psi}):
\mathcal{I}(h)=\big\{f\in\mathcal{O}_X\,;\;f\cdot\mathcal{I}(h)\subset
\mathcal{I}(he^{-\psi})\big\},
\end{equation}
and (\ref{main-surj}) can thus also be considered as a restriction morphism.
\vskip4pt plus 1pt\noindent
{\bf (B)} A surjectivity statement similar to (\ref{main-degree-q})
holds true when $(E,h)$ is a holomorphic vector bundle equipped with a smooth
hermitian \hbox{metric~$h$}. In that case, the required curvature condition
(\ref{weak-curv-cond}) is a semipositivity assumption 
\begin{equation}\label{vector-bundle-curv-cond}
\Theta_{E,h}+(1+\alpha\delta)i\ddbar \psi \otimes
\mathop{\rm Id}\nolimits_E\geq 0\qquad\text{in the sense of Nakano, for}~~
\alpha= 0,1.
\end{equation}
(This means that the corresponding hermitian form on $T_X\otimes E$ takes
nonnegative values on all tensors of $T_X\otimes E$, even those
that are non decomposable.)
\vskip4pt plus 1pt\noindent
{\bf (C)} The strength of our statements lies in the fact that no strict
positivity assumption is made. This is a typical situation in algebraic
geometry, e.g.\ in the study of the minimal model program (MMP)
for varieties which are not of general type. Our joint work [DHP13]
contains some algebraic applications which we intend
to reinvestigate in future work, by means of the present stronger
qualitative statements.
\vskip4pt plus 1pt\noindent
{\bf (D)} Notice that if one replaces
(\ref{weak-curv-cond}) by a strict positivity hypothesis
\begin{equation}\label{strict-curv-cond}
\Theta_{E,h}+i\ddbar \psi \geq \varepsilon\omega\qquad\text{in the sense
of currents, for some}~~\varepsilon>0,
\end{equation}
then Nadel's vanishing theorem implies
$H^q(X,\mathcal{O}_X(K_X\otimes E)\otimes\mathcal{I}(he^{-\psi}))=0$
for $q\geq 1$, and the injectivity and surjectivity statements are just
trivial consequences.
\vskip4pt plus 1pt\noindent
{\bf (E)} By applying convex combinations, one sees that
condition (\ref{weak-curv-cond}) takes an equivalent form if we assume
the inequality to hold for $\alpha$ varying in the whole
interval~$[0,1]$.\qed}
\end{remark}

We now turn ourselves to the problem of establishing $L^2$ estimates for
the extension problem, along the lines of \cite{OT87}. The reader will find
all details in \cite{Dem15b}.

\begin{definition} If $\psi$ is a quasi-psh function on a
complex manifold $X$, we
say that the singularities of $\psi$ are {\rm log canonical} along the
zero variety $Y=V(\mathcal{I}(e^{-\psi}))$ if~
$\mathcal{I}(e^{-(1-\varepsilon)\psi})_{|Y}
=\mathcal{O}_{X\restriction Y}$ for every $\varepsilon>0$.
\end{definition}

In case $\psi$ has log canonical singularities, it is easy to see by the
H\"older inequality and the result of Guan-Zhou \cite{GZ15b} on the
``strong openness conjecture'' that $\mathcal{I}(\psi)$ is a reduced
ideal, i.e.\ that $Y=V(\mathcal{I}(\psi))$ is a reduced
analytic subvariety of $X$. If $\omega$ is a K\"ahler
metric on~$X$, we let $dV_{X,\omega}={1\over n!}\omega^n$ be the corresponding
K\"ahler volume element, $n=\dim X$. In case $\psi$ has log canonical
singularities along $Y=V(\mathcal{I}(\psi))$, one can also associate in a
natural way a measure $dV_{Y^\circ,\omega}[\psi]$ on the set
$Y^\circ=Y_\reg$ of regular points of $Y$ as follows. If
$g\in\mathcal{C}_c(Y^\circ)$ is a compactly supported continuous function on
$Y^\circ$ and $\widetilde g$ a compactly supported extension of $g$ to
$X$, we set
\begin{equation}\label{residue-measure}
\int_{Y^\circ}g\,dV_{Y^\circ,\omega}[\psi]=
\limsup_{t\to-\infty}\int_{\{x\in X\,,\;t<\psi(x)<t+1\}}
\widetilde ge^{-\psi}\,dV_{X,\omega}.
\end{equation}
By the Hironaka desingularization theorem, one can show that the
limit does not depend on the
continuous extension~$\widetilde g$, and that one gets in this way a
measure with smooth positive density with respect to the Lebesgue
measure, at least on an (analytic) Zariski open set in $Y^\circ$.
In case $Y$ is a codimension $r$ subvariety of $X$
defined by an equation $\sigma(x)=0$ associated with a section $\sigma
\in H^0(X,S)$ of some hermitian vector bundle $(S,h_S)$ on $X$, and
assuming that $\sigma$ is generically transverse to zero along $Y$, it is
natural to take
\begin{equation}\label{psi-from-section}
\psi(z)=r\log|\sigma(z)|^2_{h_S}.
\end{equation}  
One can then easily check that $dV_{Y^\circ,\omega}[\psi]$ is the
measure supported on $Y^\circ=Y_{\rm reg}$ such~that
\begin{equation}\label{dV-from-section}
dV_{Y^\circ,\omega}[\psi]={2^{r+1}\pi^r\over (r-1)!}\,
{1\over|\Lambda^r(d\sigma)|^2_{\omega,h_S}}dV_{Y,\omega}\quad
\hbox{where}\quad
dV_{Y,\omega}= {1\over (n-r)!}\,\omega^{n-r}_{|Y^\circ}.
\end{equation}
For a quasi-psh function with log canonical singularities, 
$dV_{Y^\circ,\omega}[\psi]$ should thus be seen as some sort of
(inverse of) Jacobian determinant associated with the logarithmic 
singularities of~$\psi$. In general, the measure $dV_{Y^\circ,\omega}[\psi]$
blows up (i.e.\ has infinite volume) in a neighborhood of
singular points of $Y$. Finally, the following positive real function
will make an appearance in several of our estimates~:
\begin{equation}
\gamma(x) = \exp(-x/2)~~\hbox{if $x\ge 0$},\qquad
\gamma(x) = \frac{1}{1+x^2}~~\hbox{if $x\le 0$}.
\end{equation}
The first generalized $L^2$ estimate  we are interested in
is a variation of Theorem~4 in [Ohs01]. One difference
is that we do not require any specific behavior of the quasi-psh function
$\psi$ defining the subvariety: any quasi-psh function with log canonical 
singularities will do; secondly, we do not want to make any assumption
that there exist negligible sets in the ambient manifold whose complements 
are Stein, because such an hypothesis need not be true
on a general compact K\"ahler manifold -- one of the targets of our study.

\begin{thm}[$L^2$ estimate for the extension from reduced subvarieties]
\label{reduced-estimate}
Let $X$ be a holomorphically convex K\"ahler manifold, and $\omega$ a
K\"ahler metric on~$X$.  Let $(E,h)$ be a holomorphic vector bundle 
equipped with a smooth hermitian metric $h$ on $X$, and let
$\psi:X\to[-\infty,+\infty[$ be a quasi-psh function on $X$ with neat
analytic singularities. Let $Y$ be the analytic subvariety of $X$ defined 
by $Y=V(\mathcal{I}(e^{-\psi}))$ and assume that $\psi$ has log canonical
singularities along $Y$, so that $Y$ is reduced. Finally, assume
that the Chern curvature tensor $\Theta_{E,h}$ is such that the sum
$$
\Theta_{E,h}+(1+\alpha\delta)\,i\ddbar\psi\otimes\mathop{\rm Id}\nolimits_E
$$
is Nakano semipositive for
some $\delta>0$ and $\alpha=0,1$. Then for every holomorphic
section $f\in H^0(Y^\circ,(K_X\otimes E)_{|Y^\circ})$ on $Y^\circ=Y_\reg$
such that
$$
\int_{Y^\circ}|f|^2_{\omega,h}dV_{Y^\circ,\omega}[\psi]<+\infty,
$$
there exists an extension $F\in H^0(X,K_X\otimes E)$ whose restriction to
$Y^\circ$ is equal to $f$, such that
$$
\int_X\gamma(\delta\psi)\,|F|^2_{\omega,h}e^{-\psi}
dV_{X,\omega}\le{34\over\delta}\int_{Y^\circ}|f|^2_{\omega,h}dV_{Y^\circ,\omega}[\psi].
$$
\end{thm}

\begin{remark} {\rm Although $|F|^2_{\omega,h}$ and
$dV_{X,\omega}$ both depend on $\omega$, it is easy to see that the
product $|F|^2_{\omega,h}dV_{X,\omega}$ actually does not depend on
$\omega$ when $F$ is a $(n,0)$-form. The same observation applies to
the product $|f|^2_{\omega,h}dV_{Y^\circ,\omega}[\psi]$, hence the
final $L^2$ estimate is in fact independent of $\omega$. Nevertheless,
the existence of a K\"ahler metric (and even of a complete K\"ahler
metric) is crucial in the proof, thanks to the techniques developped in
[AV65] and [Dem82]. The constant 34 is of course non optimal;
the technique developed in \cite{GZ15a} provides optimal
choices of the function $\gamma$ and of the constant in the right hand side.
\qed}
\end{remark}

We now turn ourselves to the case where non reduced multiplier ideal sheaves
and non reduced subvarieties are considered. This situation has already
been considered by D.~Popovici [Pop05] in the case of powers of a
reduced ideal, but we aim here at a much wider generality, which also
yields more natural assumptions. For $m\in\mathbb{R}_+$, we consider the
multiplier ideal sheaf $\mathcal{I}(e^{-m\psi})$ and the associated
non necessarily reduced
subvariety $Y^{(m)}=V(\mathcal{I}(e^{-m\psi}))$, together with the
structure sheaf
$\mathcal{O}_{Y^{(m)}}=\mathcal{O}_X/\mathcal{I}(e^{-m\psi})$,
the real number $m$ being viewed as some sort
of multiplicity -- the support $|Y^{(m)}|$ may increase with $m$, but certainly
stabilizes to the set of poles $P=\psi^{-1}(-\infty)$ for $m$ large enough. 
We assume the existence of a discrete sequence of positive numbers
$$0=m_0<m_1<m_2<\ldots<m_p<\ldots$$
such that $\mathcal{I}(e^{-m\psi})=\mathcal{I}(
e^{-m_p\psi})$ for $m\in [m_p,m_{p+1}[$ (with of course
$\mathcal{I}(e^{-m_0\psi})=\mathcal{O}_X)$; they
are called the {\it jumping numbers} of~$\psi$. The existence of a discrete
sequence of jumping numbers is automatic if $X$ is compact. In general, this
still holds on every relatively compact open subset 
$$X_c:=\{x\in X\,,\;\rho(x)<c\}\Subset X,$$ 
but requires some some of uniform behaviour of singularities at infinity
in the non compact case. We are interested in extending a holomorphic section 
$$
f\in H^0(Y^{(m_p)},\mathcal{O}_{Y^{(m_p)}}(K_X\otimes E_{|Y^{(m_p)}})
:=H^0(Y^{(m_p)},\mathcal{O}_X(K_X\otimes_{\mathbb{C}}E)
\otimes_{\mathcal{O}_X}\mathcal{O}_X/\mathcal{I}(e^{-m_p\psi})).
$$
[Later on, we usually omit to specify the rings over which tensor products
are taken, as they are implicit from the nature of objects under consideration].
The results are easier to state in case one takes a nilpotent section of
the form
$$
f\in H^0(Y^{(m_p)},\mathcal{O}_X(K_X\otimes E)\otimes\mathcal{I}(e^{-m_{p-1}\psi})/\mathcal{I}(e^{-m_p\psi})).
$$
Then $\mathcal{I}(e^{-m_{p-1}\psi})/\mathcal{I}(e^{-m_p\psi}))$ is
actually a coherent sheaf, and one can see that its
support is a reduced subvariety $Z_p$ of $Y^{(m_p)}$. Therefore 
$\mathcal{I}(e^{-m_{p-1}\psi})/\mathcal{I}(e^{-m_p\psi}))$ can be seen
as a vector bundle
over a Zariski open set $Z_p^\circ\subset Z_p$. We can mimic
formula (\ref{residue-measure})
and define some sort of infinitesimal ``$m_p$-jet'' $L^2$ norm
$|J^{m_p}f|^2_{\omega,h}\,dV_{Z_p^\circ,\omega}[\psi]$ (a purely formal
notation), as the measure on $Z_p^\circ$ defined by
\begin{equation}
\int_{Z_p^\circ}g\,|J^{m_p}f|^2_{\omega,h}\,dV_{Z_p^\circ,\omega}[\psi]=
\limsup_{t\to-\infty}\int_{\{x\in X\,,\;t<\psi(x)<t+1\}}
\widetilde g\,|\widetilde f|^2_{\omega,h}e^{-m_p\psi}\,dV_{X,\omega}
\end{equation}
for any $g\in\mathcal{C}_c(Z_p^\circ)$, where $\widetilde g\in\mathcal{C}_c(X)$ is a 
continuous extension of $g$ and $\widetilde f$ a smooth
extension of $f$ on $X$ such that $\smash{\widetilde f} - f\in
\mathcal{I}(m_p\psi)\otimes_{\mathcal{O}_X}\mathcal{C}^\infty$
(this measure again has a
smooth positive density on a Zariski open set in $Z_p^\circ$, and does
not depend on the choices of $\smash{\widetilde f}$ and~$\widetilde g$).
We extend the measure as being $0$ on 
$\smash{Y^{(m_p)}_{\rm red}}\smallsetminus Z_p$,
since $\mathcal{I}(e^{-m_{p-1}\psi})/\mathcal{I}(e^{-m_p\psi}))$ 
has support
in~$Z_p^\circ\subset Z_p$. In this context, we introduce the following
natural definition.

\begin{definition} We define the {\rm restricted multiplied ideal sheaf}
$$
\mathcal{I}'(e^{-m_{p-1}\psi})\subset\mathcal{I}(e^{-m_{p-1}\psi})
$$
to be the set of germs $F\in\mathcal{I}(e^{-m_{p-1}\psi})_x\subset
\mathcal{O}_{X,x}$ such that 
there exists a neighborhood $U$ of $x$ satisfying
$$
\int_{Y^{(m_p)}\cap U}|J^{m_p}F|^2_{\omega,h}\,
dV_{Y^{(m_p)},\omega}[\psi]<+\infty.
$$
This is a coherent ideal sheaf that contains $\mathcal{I}(e^{-m_p\psi})$.
Both of the inclusions
$$
\mathcal{I}(e^{-m_p\psi})\subset \mathcal{I}'(e^{-m_{p-1}\psi})\subset
\mathcal{I}(e^{-m_{p-1}\psi})
$$
can be strict $($even for $p=1)$.
\end{definition}

One of the geometric consequences is the following ``quantitative''
surjectivity statement, which is the analogue of
Theorem~$\ref{reduced-estimate}$ for the case when the first non
trivial jumping number $m_1=1$ is replaced by a higher jumping
number~$m_p$.

\begin{thm}\label{estimate-for-last-jump}
With the above notation and in the general 
setting of Theorem~$\ref{reduced-estimate}$ $($but without the hypothesis 
that the quasi-psh function $\psi$ has log canonical singularities$)$, let
$0=m_0<m_1<m_2<\ldots<m_p<\ldots$ be the jumping numbers of~$\psi$. Assume 
that
$$
\Theta_{E,h}+i(m_p+\alpha\delta)
\ddbar\psi\otimes\mathop{\rm Id}\nolimits_E\ge 0
$$
is Nakano semipositive for $\alpha=0,1$ and some $\delta>0$.
\begin{itemize}
\item[(a)] Let
$$
f\in H^0(Y^{(m_p)},\mathcal{O}_X(K_X\otimes E)\otimes\mathcal{I}'
(e^{-m_{p-1}\psi})/\mathcal{I}(e^{-m_p\psi}))
$$
be a section such that
$$
\int_{Y^{(m_p)}}|J^{m_p}f|^2_{\omega,h}\,dV_{Y^{(m_p)},\omega}[\psi]<+\infty.
$$
Then there exists a global section
$$
F\in H^0(X,\mathcal{O}_X(K_X\otimes E)\otimes\mathcal{I}'(e^{-m_{p-1}\psi}))
$$
which maps to $f$ under the morphism $\mathcal{I}'(e^{-m_{p-1}\psi})\to
\mathcal{I}(e^{-m_{p-1}\psi})/\mathcal{I}(e^{-m_p\psi})$, such that
$$
\int_{X}\gamma(\delta\psi)\,|F|^2_{\omega,h}\,e^{-m_p\psi}dV_{X,\omega}[\psi]
\le {34\over\delta}\int_{Y^{(m_p)}}|J^{m_p}f|^2_{\omega,h}\,dV_{Y^{(m_p)},\omega}[\psi].
$$
\item[(b)] The restriction morphism
$$
\begin{matrix}  
&\kern-100pt
H^0(X,\mathcal{O}_X(K_X\otimes E)\otimes\mathcal{I}'(e^{-m_{p-1}\psi}))\\
\noalign{\vskip5pt}
&\kern100pt{}
\to  H^0(Y^{(m_p)},\mathcal{O}_X(K_X\otimes E)\otimes\mathcal{I}'(e^{-m_{p-1}\psi})/\mathcal{I}(e^{-m_p\psi}))
\end{matrix}
$$
is surjective.
\end{itemize}
\end{thm}

If $E$ is a line bundle and $h$ a singular hermitian metric on $E$, a
similar result can be obtained by approximating $h$. However, the
$L^2$ estimates then require to incorporate $h$ into the
definition of the multiplier ideals, as in
Theorem~\ref{main-theorem} (see \cite{Dem15a}). Hosono \cite{Hos17a} has
shown that one can obtain again an optimal $L^2$ estimate in the situation
of Theorem~\ref{estimate-for-last-jump}, when $\mathcal{I}(e^{-m_p\psi})$ is
  a power of the reduced ideal of $Y$.

\begin{question} {\rm It would be interesting to know whether
Theorem~\ref{main-theorem} can be strengthened by suitable
$L^2$ estimates, without making undue additional hypotheses
on the section $f$ to extend. The main difficulty is already to define the 
norm of jets when there is more than one jump number involved. Some sort of
``Cauchy inequality'' for jets would be needed in order to derive the
successive jet norms from a known global $L^2$ estimate for a holomorphic
section defined on the whole of~$X$. We do not know how to proceed
further at this point.}
\end{question}

\vskip8pt

\section{Bochner-Kodaira estimate with approximation}

The crucial idea of the proof is to prove the results (say, in the form of the surjectivity statement), only up to approximation. This is done by solving a $\dbar$-equation
\begin{equation*}
\dbar u_\varepsilon+w_\varepsilon=v
\end{equation*}
where the right hand side $v$ is given and $w_\varepsilon$ is an error term such
that $\Vert w_\varepsilon\Vert=O(\varepsilon^a)$ as $\varepsilon\to 0$,
for some constant $a>0$. A twisted Bochner-Kodaira-Nakano identity introduced by
Donnelly and Fefferman \cite{DF83}, and Ohsawa and Takegoshi \cite{OT87} is used for that purpose. The technology goes back
to the fundamental work of Bochner (\cite{Boc48}), Kodaira
(\cite{Kod53a,Kod53b,Kod54}), Akizuki-Nakano (\cite{AN54}, \cite{Nak55}),
Kohn (\cite{FK}), Andreotti-Vesentini (\cite{AV65}),
H\"ormander (\cite{Hor65,Hor66}). The version 
we need uses in an essential way an additional correction term, so as
to allow a weak positivity hypothesis. It can be stated as follows.

\begin{proposition}{\rm (see \cite[Prop.~3.12]{Dem15b})}\label{L2-estimate}
Let $X$ be a complete K\"ahler manifold
equipped with a $($non necessarily complete$)$ K\"ahler metric
$\omega$, and let $(E,h)$ be a Hermitian vector bundle over~$X$.
Assume that there are smooth and bounded functions $\eta,\,\lambda>0$
on $X$ such that the curvature operator
$$
B=B^{n,q}_{E,h,\omega,\eta,\lambda}=
[\eta\,\Theta_{E,h}-i\,\ddbar\eta-i\lambda^{-1}d\partial\eta\wedge\dbar\eta,
\Lambda_\omega]\in C^\infty(X,\mathop{\rm Herm}(\Lambda^{n,q}T^*_X\otimes E))
$$
satisfies $B+\varepsilon I>0$ for some $\varepsilon>0$ $($so that $B$
can be just semi-positive or even slightly negative; here $I$ is the identity 
endomorphism$)$. Given a section
$v\in L^2(X,\Lambda^{n,q}T^*_X\otimes E)$ such that $\dbar v=0$ and
$$M(\varepsilon):=
\int_X\langle (B+\varepsilon I)^{-1}v,v\rangle\,dV_{X,\omega}<+\infty,$$
there exists an approximate solution
\hbox{$f_\varepsilon\in L^2(X,\Lambda^{n,q-1}T^*_X \otimes E)$} and a correction
term $w_\varepsilon\in L^2(X,\Lambda^{n,q}T^*_X \otimes E)$ such that
$\dbar u_\varepsilon=v-w_\varepsilon$ and
$$
\int_X(\eta+\lambda)^{-1}|u_\varepsilon|^2\,dV_{X,\omega}+
\frac{1}{\varepsilon}\int_X|w_\varepsilon|^2\,dV_{X,\omega}\le M(\varepsilon).
$$
Moreover, if $v$ is smooth, then $u_\varepsilon$ and $w_\varepsilon$ can be taken smooth.
\end{proposition}

In our situation, the main part of the solution, namely $u_\varepsilon$, may very well explode as $\varepsilon\to 0$. In order to show that the equation $\dbar u=v$ can be solved, it is therefore needed to check that the space of coboundaries is closed in the space of cocycles in the Fr\'echet topology under consideration (here, the $L^2_{\rm loc}$ topology), in other words, that the related cohomology group
$H^q(X,\mathcal{F})$ is Hausdorff. In this respect, the fact of considering
$\dbar$-cohomology of smooth forms equipped with the $C^\infty$ topology
on the one hand, or cohomology of forms $u\in L^2_{\rm loc}$ with
$\dbar u\in L^2_{\rm loc}$ on the other hand, yields the same topology on
the resulting  cohomology group $H^q(X,\mathcal{F})$. This comes from the fact
that both complexes yield fine resolutions of the same coherent sheaf
$\mathcal{F}$, and the topology of $H^q(X,\mathcal{F})$ can also be obtained
by using \v{C}ech cochains with respect  to a Stein covering $\mathcal{U}$
of~$X$. The required Hausdorff property then comes from the following well
known fact.

\begin{lemma}\label{hausdorff-lemma}Let $X$ be a holomorphically convex
complex space and $\mathcal{F}$ a coherent analytic sheaf over $X$.
Then all cohomology groups $H^q(X,\mathcal{F})$ are Hausdorff with respect
to their natural topology $($induced by the Fr\'echet topology of local uniform
convergence of holomorphic 
cochains$)$.
\end{lemma}

In fact, the Remmert reduction theorem implies that $X$ admits a proper
holomorphic map $\pi:X\to S$ onto a Stein space $S$, and Grauert's
direct image theorem shows that all direct images $R^q\pi_*\mathcal{F}$
are coherent sheaves on~$S$. Now, as $S$ is Stein, Leray's theorem combined
with Cartan's theorem B tells us that we have
an isomorphism $H^q(X,\mathcal{F})\simeq H^0(S,R^q\pi_*\mathcal{F})$. More
generally, if $U\subset S$ is a Stein open subset, we have
\begin{equation}\label{local-over-S}
H^q(\pi^{-1}(U),\mathcal{F})\simeq H^0(U,R^q\pi_*\mathcal{F})
\end{equation}
and when $U\Subset S$ is relatively compact,
it is easily seen that this a topological isomorphism of Fr\'echet spaces
since both sides are $\mathcal{O}_S(U)$ modules of finite type and can be
seen as a Fr\'echet quotient of some direct sum
$\mathcal{O}_S(U)^{\oplus N}$ by looking at local generators and local relations
of~$R^q\pi_*\mathcal{F}$.
Therefore $H^q(X,\mathcal{F})\simeq H^0(S,R^q\pi_*\mathcal{F})$
is a topological isomorphism and the space of sections in the right
hand side is a Fr\'echet space. In particular,
$H^q(X,\mathcal{F})$ is Hausdorff.\qed

\section{Sketch of proof of the extension theorem}

The reader may consult \cite{Dem15b} and \cite{CDM17} for more details.
After possibly shrinking $X$ into a relatively compact holomorphically
convex open subset $X'=\pi^{-1}(S')\Subset X$, we can suppose that
$\delta>0$ is a constant and that $\psi\leq 0$ (otherwise subtract
a large constant to~$\psi$). As $\pi:X\to S$ is proper, we can also
assume that $X$ admits a finite Stein covering $\mathcal{U}=(U_i)$.
Any cohomology class in
$$
H^q(Y,\mathcal{O}_X(K_X\otimes E)\otimes\mathcal{I}(h)/\mathcal{I}(he^{-\psi}))
$$
is represented by a holomorphic \v{C}ech $q$-cocycle with respect to
the covering $\mathcal{U}$
$$
(c_{i_0\ldots i_q}),\qquad
c_{i_0\ldots i_q}\in H^0\big(U_{i_0}\cap\ldots\cap U_{i_q},
\mathcal{O}_X(K_X\otimes E)\otimes\mathcal{I}(h)/\mathcal{I}(he^{-\psi})\big).
$$
By the standard sheaf theoretic isomorphisms with Dolbeault cohomology (cf.\
e.g.\ \cite{Dem-e-book}), this 
class is represented by a smooth $(n,q)$-form
$$
f=\sum_{i_0,\ldots,i_q}c_{i_0\ldots i_q}\rho_{i_0}
\dbar\rho_{i_1}\wedge\ldots\dbar\rho_{i_q}
$$
by means of a partition of unity $(\rho_i)$ subordinate to $(U_i)$. This form is
to be interpreted as a form on the (non reduced) analytic subvariety $Y$ associated with the ideal sheaf $\mathcal{J}=\mathcal{I}(he^{-\psi}):
\mathcal{I}(h)$ and the structure sheaf $\mathcal{O}_Y=
\mathcal{O}_X/\mathcal{J}$. We get an extension as a smooth (no longer
$\dbar$-closed) $(n,q)$-form on $X$ by taking
$$
\widetilde{f}=\sum_{i_0,\ldots,i_q}\widetilde{c}_{i_0\ldots i_q}\rho_{i_0}
\dbar\rho_{i_1}\wedge\ldots\dbar\rho_{i_q}
$$
where $\widetilde{c}_{i_0\ldots i_q}$ is an extension of
$c_{i_0\ldots i_q}$ from $U_{i_0}\cap\ldots\cap U_{i_q}\cap Y$ to
$U_{i_0}\cap\ldots\cap U_{i_q}$.
Without loss of generality, we can assume that $\psi$ admits 
a discrete sequence of ``jumping numbers''
\begin{equation}\label{jumping-numbers}
0=m_0<m_1<\cdots <m_p<\cdots\quad\hbox{
such that $\mathcal{I}(m\psi)=\mathcal{I}(m_p\psi)$ for
$m\in[m_p,m_{p+1}[$}.
\end{equation}
 Since $\psi$ is assumed to have analytic singularities, 
this follows from using a log resolution of singularities, thanks to the
Hironaka desingularization theorem (by the much deeper result of
\cite{GZ15b} on the strong openness conjecture, one could even possibly 
eliminate the assumption that $\psi$ has analytic singularities).
We fix here $p$ such that $m_p\leq 1<m_{p+1}$, and
in the notation of \cite{Dem15b}, we let $Y=Y^{(m_p)}$ be defined by the
non necessarily reduced
structure sheaf $\mathcal{O}_Y=\mathcal{O}_X/\mathcal{I} (e^{-\psi})=
\mathcal{O}_X/\mathcal{I} (e^{-m_p\psi})$.

We now explain the choice of metrics and auxiliary functions $\eta$, $\lambda$
for the application of Proposition~\ref{L2-estimate},
following the arguments of \cite[proof of th.~2.14, p.~217]{Dem15b}.
Let $t\in \mathbb{R}^-$ and let $\chi_t$ be the negative convex increasing function defined in \cite[(5.8$\,*$), p.~211]{Dem15b}. Put $\eta_t := 1 -\delta \cdot \chi_t(\psi)$
and $\lambda_t := 2 \delta \frac{(\chi_t^2(\psi))^2}{\chi_t''(\psi)}$. We set
\begin{eqnarray*}
R_t &:=& \eta_t (\Theta_{E,h} +i\ddbar\psi ) -i\ddbar \eta_t -\lambda_t ^{-1} i \partial \eta_t \wedge \dbar \eta_t\\
&\kern3pt =& \eta_t (\Theta_{E,h} + (1+\delta \eta_t ^{-1} \chi_t'(\psi))i\ddbar\psi ) 
+ \frac{\delta\cdot \chi_t''(\psi)}{2}  i \partial\psi\wedge \dbar\psi.
\end{eqnarray*}
Note that $\chi_t''(\psi) \geq \frac{1}{8}$ on $W_t =\{t < \psi < t+1\}$. The curvature assumption \eqref{weak-curv-cond} implies
$$\Theta_{E,h} + (1+\delta \eta_t ^{-1} \chi_t'(\psi))\,i\ddbar\psi \geq 0 \qquad\text{on }X .$$
As in \cite{Dem15b}, we find
\begin{equation}\label{posit-curv-bound1}
R_t \geq 0 \qquad\text{on~~}X
\end{equation}
and
\begin{equation}\label{posit-curv-bound2}
R_t \geq \frac{\delta}{16} i \partial\psi\wedge \dbar\psi\qquad\text{on~~}W_t =\{t < \psi < t+1\}. 
\end{equation}
\vskip6pt plus 1pt minus 1pt

\noindent
Let $\theta : [ -\infty , +\infty [{}\rightarrow [0,1]$ be a smooth non increasing real function satisfying
$\theta (x)=1$ for $x \leq 0$, $\theta (x)=0$ for $x \geq 1$ and $|\theta' | \leq 2$. By using a blowing up process, one can reduce the situation to the case where $\psi$ has divisorial singularities. Then  we still have
$$\Theta_{E,h}+ (1+\delta \eta_t ^{-1} \chi_t'(\psi)) 
(i\ddbar\psi)_{\ac} \geq 0 \qquad\text{on }X ,$$
where $(i\ddbar\psi)_{\ac}$ is the absolutely continuous part of
$i\ddbar\psi$.  The regularization techniques of \cite{DPS01} and
\cite[Th.~1.7, Remark~1.11]{Dem15a} produce a family of singular
metrics $\{h_{t,\varepsilon}\}_{k=1}^{+\infty}$ which are smooth in
the complement $X\smallsetminus Z_{t,\varepsilon}$ of an analytic set,
such that $\mathcal{I} (h_{t,\varepsilon}) =\mathcal{I} (h)$,
$\mathcal{I} (h_{t,\varepsilon}e^{-\psi}) =\mathcal{I} (h e^{-\psi})$ and
$$
\Theta_{E,h_{t,\varepsilon}}+
(1+\delta \eta_t ^{-1} \chi_t'(\psi))\,i\ddbar\psi \geq
-\frac{1}{2}\varepsilon \omega \qquad\text{on }X .
$$
The additional error term $-\frac{1}{2}\varepsilon \omega$ is
irrelevant when we use Proposition \ref{L2-estimate}, as it is absorbed by
taking the hermitian operator $B+\varepsilon I$. Therefore for every
$t\in \mathbb{R}^-$, with the adjustment $\varepsilon=e^{\alpha t}$,
$\alpha\in{}]0,m_{p+1}-1[$, we can find a singular metric
$h_t=h_{t,\varepsilon}$ which is smooth in the complement $X\setminus Z_t$
of an analytic set, such that
$\mathcal{I} (h_t) =\mathcal{I} (h)$,
$\mathcal{I} (h_te^{-\psi}) =\mathcal{I} (he^{-\psi})$
and $h_t\uparrow h$ as
$t\rightarrow -\infty$. We now apply the $L^2$ estimate of 
Proposition \ref{L2-estimate} and observe that $X\smallsetminus Z_t$ is complete K\"ahler (at least after we shrink $X$ a little bit as $X'=\pi^{-1}(S')$,
cf.~\cite{Dem82}). As a consequence, one can find  sections 
$u_t$, $w_t$ satisfying 
\begin{equation}\label{l2estimate}
\dbar u_t+ w_t=v_t:=
\dbar\big(\theta(\psi -t) \cdot \widetilde{f}\;\big)
\end{equation}
and
\begin{equation}\label{main-l2-estimate}
\begin{matrix}
&\displaystyle
  \int_X (\eta_t+\lambda_t)^{-1}| u_t|_{\omega, h_t} ^2 e^{-\psi} dV_{X,\omega}+
\frac{1}{\varepsilon} \int_X |w_t|_{\omega, h_t}^2 e^{-\psi} 
dV_{X,\omega}\\
&\displaystyle{}\kern40pt{}\leq 
\int_X \langle (R_t +\varepsilon I)^{-1}v_t,v_t\rangle_{\omega, h_t}\;
e^{-\psi}dV_{X,\omega}.
\end{matrix}
\end{equation}
One of the main consequence of (\ref{posit-curv-bound2}) and
(\ref{main-l2-estimate}) is that, for
$\varepsilon=e^{\alpha t}$ and $\alpha$ well chosen, one can infer that
the error term satisfies
$$ \lim_{t\rightarrow -\infty} \int_X |w_t|_{\omega, h_t} ^2 e^{-\psi} d V_{X,\omega} =0 .$$
One difficulty, however, is that $L^2$ sections cannot be restricted in a
continuous
way to a subvariety. In order to overcome this problem,
we play again the game of returning to \v{C}ech cohomology by solving 
inductively $\dbar$-equations for $w_t$ on $U_{i_0}\cap\ldots\cap U_{i_k}$,
until we reach an equality
\begin{equation}\label{l2-estimate-nq-forms}
\dbar\big(\theta(\psi -t) \cdot \widetilde{f}  -\widetilde{u}_t\big)
=\widetilde{w}_t:=
-\sum_{i_0,\ldots,i_{q-1}}s_{t,i_0\ldots i_q}\dbar\rho_{i_0}\wedge
\dbar\rho_{i_1}\wedge\ldots\dbar\rho_{i_q}
\end{equation}
with holomorphic sections 
$s_{t,I}=s_{t,i_0\ldots i_q}$ on $U_I=U_{i_0}\cap\ldots\cap U_{i_q}$,
such that
$$ \lim_{t\rightarrow -\infty} \int_{U_I} |s_{t,I}|_{\omega, h_t} ^2 
e^{-\psi} d V_{X,\omega} =0.$$
Then the right hand side of (\ref{l2-estimate-nq-forms}) is smooth, and more
precisely has coefficients in the sheaf
$\mathcal{C}^\infty\otimes_{\mathcal{O}}\mathcal{I}(he^{-\psi})$, and 
$\widetilde{w}_t\to 0$ in $C^\infty$ topology. A~priori,
$\widetilde{u}_t$ is an $L^2$ $(n,q)$-form equal to $u_t$ plus a combination
$\sum\rho_is_{t,i}$ of the local solutions of $\dbar s_{t,i}=w_t$, plus
$\sum\rho_is_{t,i,j}\wedge\dbar\rho_j$ where
$\dbar s_{t,i,j}=s_{t,j}-s_{t,i}$, plus etc~$\ldots$~, and is such that
$$\int_{X} |\widetilde{u}_t|_{\omega, h_t}^2 
e^{-\psi} d V_{X,\omega} <+\infty.$$
Since $H^q(X,\mathcal{O}_X(K_X\otimes E)\otimes\mathcal{I}(he^{-\psi}))$ can be 
computed with the $L^2_{\rm loc}$ resolution of the coherent sheaf, or 
alternatively  with the $\dbar$-complex of $(n,{\scriptstyle\bullet})$-forms 
with coefficients in
$\mathcal{C}^\infty\otimes_{\mathcal{O}}\mathcal{I}(he^{-\psi})$, we may assume that
\hbox{$\widetilde{u}_t\in\mathcal{C}^\infty\otimes_{\mathcal{O}}
\mathcal{I}(he^{-\psi})$},
after playing again with \v{C}ech cohomology. Lemma~\ref{hausdorff-lemma}
yields a sequence of smooth $(n,q)$-forms $\sigma_t$ with coefficients in 
$\mathcal{C}^\infty\otimes_{\mathcal{O}}\mathcal{I}(h)$,
such that $\dbar\sigma_t=\widetilde{w}_t$ and $\sigma_t\to 0$ in
$C^\infty$-topology. Then $\widetilde{f}_t=
\theta(\psi -t) \cdot \widetilde{f}  -\widetilde{u}_t-\sigma_t$
is a $\dbar$-closed $(n,q)$-form on $X$ with values in
$\mathcal{C}^\infty\otimes_{\mathcal{O}}\mathcal{I}(h)\otimes
\mathcal{O}_X(E)$, whose image in
$H^q(X,\mathcal{O}_X(K_X\otimes E)\otimes\mathcal{I}(h)/\mathcal{I}(he^{-\psi}))$
converges to $\{f\}$ in $C^\infty$ Fr\'echet topology. We conclude by a
density argument on the Stein space $S$, by looking at the coherent 
sheaf morphism
$$R^q\pi_*\big(\mathcal{O}_X(K_X\otimes E)\otimes\mathcal{I}(h)\big)\to
R^q\pi_*\big(\mathcal{O}_X(K_X\otimes E)\otimes\mathcal{I}(h)/
\mathcal{I}(he^{-\psi})\big).\eqno\qed
$$
\vskip4pt

\noindent
{\bf Proof of the quantitative estimates.} We refer again to \cite{Dem15b}
for details. One of the main features of the above qualitative proof 
is that we have not tried to control the solution $u_t$
of our $\dbar$-equation, in fact we only needed to prove that the error
term $w_t$ converges to zero. However, to get quantitative
$L^2$ estimates, we have to pay attention to the $L^2$ norm
of $u_t$. It is under control as
$t\to-\infty$ only when $f$ satisfies the more restrictive
condition of being $L^2$ with respect to the residue measure
$dV_{Y^\circ,\omega}[\psi]$. This is the reason why we lose track
of the solution when the volume of the measure explodes on
$Y_{\rm sing}$, or when there are several jumps involved in the
multiplier ideal sheaves.

\section{Applications of the Ohsawa-Takegoshi
extension theorem}

The Ohsawa-Takegoshi extension theorem is a very powerful tool that
has many important applications to complex analysis and geometry.
We will content ourselves by mentioning only a few statements and references.

\subsection{Approximation of plurisubharmonic functions and
  of closed $(1,1)$-currents.} By considering the extension from points
(i.e.\ a $0$-dimensional connected subvariety $Y\subset X$), even just locally
on coordinates balls, one gets a precise Bergman kernel estimate for 
Hilbert spaces attached to multiples of any plurisubharmonic function.
This leads to regularization theorems (\cite{Dem92}) that have many
applications, such as the Hard Lefschetz theorem with multiplier ideal
sheaves (\cite{DPS01}), or extended vanishing theorems for pseudoeffective
line bundles (\cite{Cao14}). The result may consult \cite{Dem15a} for a
survey of these questions.  Another consequence is a very simple and
direct proof of Siu's result \cite{Siu74} on the analyticity of
sublevel sets of Lelong numbers of closed positive currents.

\subsection{Invariance of plurigenera.} Around 2000, Siu
\cite{Siu02} proved that
for every smooth projective deformation $\pi:\mathcal{X}\to S$ over an
irreducible base~$S$, the plurigenera $p_m(t)=h^0(X_t,K_{X_t}^{\otimes m})$
of the fibers $X_t=\pi^{-1}(t)$ are constant. The proof relies in an
essential way on the Ohsawa-Takegoshi extension theorem, and was later
simplified and generalized by P\u{a}un \cite{Pau07}. It is remarkable that no
algebraic proof of this purely algebraic result is known!

\subsection{Semicontinuity of log singularity exponents.} In
\cite{DK01}, we proved that the log singularity exponent
(or log canonical threshold) $c_x(\varphi)$, defined as
the supremum of constants $c>0$ such that
$e^{-c\varphi}$ is integrable in a neighborhood of a point~$x$, is
a lower semicontinuous function with respect to the topology
of weak convergence on plurisubharmonic functions. Guan and Zhou
\cite{GZ15b} recently proved our ``strong openness conjecture'',
namely that the integrability of $e^{-\varphi}$ implies the
integrability of $e^{-(1+\varepsilon)\varphi}$ for $\varepsilon>0$~small;
later alternative proofs have been exposed in \cite{Hie14} and \cite{Lem14}.

\subsection{Proof of the Suita conjecture.} In \cite{Blo13} B{\l}ocki
determined the value of the optimal constant in the Ohsawa-Takegoshi
extension theorem, a result that was subsequently generalized
by Guan and Zhou \cite{GZ15a}. In complex dimension~1, this result implies
in its turn a conjecture of N.~Suita, stating that for any bounded domain
$D$ in $\mathbb{C}$, one has  $c_D^2\leq\pi K_D$, where $c_D(z)$ is the
logarithmic capacity of $\mathbb{C}\smallsetminus D$ with respect to
$z\in D$ and $K_D$ is the Bergman kernel on the diagonal. Guan and Zhou
\cite{GZ15a} proved that the equality occurs if ond only if $D$ is 
conformally equivalent to the disc minus a closed set of inner capacity
zero.
\vskip18pt


\bigskip

\begin{thebibliography}{n}

\bibitem[AN54]{AN54} Akizuki, Y., Nakano, S., \textit{Note on Kodaira-Spencer's
proof of Lefschetz theorems,} Proc.\ Jap.\ Acad.\ {\bf 30} (1954) 266--272.

\bibitem[AV65]{AV65} Andreotti, A., Vesentini, E., \textit{Carleman estimates for the
Laplace-Beltrami equation in complex manifolds,} Publ.\ Math.\ I.H.E.S.\
{\bf 25} (1965) 81--130.

\bibitem[Ber96]{Ber96} Berndtsson, B., \textit{The extension theorem of Ohsawa-Takegoshi
and the theorem of Donnelly-Fefferman,} Ann.\ Inst.\ Fourier 
{\bf 14} (1996) 1087--1099.

\bibitem[BL16]{BL16} Berndtsson, B., Lempert, L., \textit{A proof of the 
Ohsawa-Takegoshi theorem with sharp estimates,}
J.\ Math.\ Soc.\ Japan {\bf 68} (2016), no.~4, 1461--1472.

\bibitem[Blo13]{Blo13} B{\l}ocki, Z., \textit{Suita conjecture and the
Ohsawa-Takegoshi extension theorem,} Invent.\ Math.\ {\bf 193} (2013),
no.~1, 149--158.

\bibitem[Boc48]{Boc48} Bochner, S.,
\textit{Curvature and Betti numbers (I) and (II)},
Ann.\ of Math.\ {\bf 49} (1948), 379--390$\,$; {\bf 50} (1949), 77--93.

\bibitem[Cao14]{Cao14} Cao, Junyan,
\textit{Numerical dimension and a Kawamata–Viehweg–Nadel-type vanishing theorem on compact K\"ahler manifolds,} Compositio Mathematica {\bf 150}
(2014) 1869--1902.
  
\bibitem[CDM17]{CDM17} Cao, Junyan, Demailly, J.-P.,
Matsumura, Shin-ichi, \textit{A general extension theorem for cohomology classes on non reduced analytic subspaces,}
Science China Mathematics,
Volume 60, Issue 6, June 2017, 949--962.

\bibitem[Che11]{Che11} Chen, Bo-Yong, \textit{A simple proof of
the Ohsawa-Takegoshi extension theorem,} arXiv: math.CV/ 1105.2430.

\bibitem[Dem82]{Dem82}
Demailly, J.-P.,
\textit{Estimations $L^{2}$ pour l'op\'erateur $\overline{\partial}$ d'un 
fibr\'e vectoriel holomorphe semi-positif au-dessus d'une vari\'et\'e k\"ahl\'erienne compl\`ete, }
Ann. Sci. \'Ecole Norm. Sup(4). {\bf{15}} (1982), 457--511. 

\bibitem[Dem92]{Dem92} Demailly, J.-P., \textit{Regularization of closed
positive currents and Intersection Theory,}
J.\ Alg.\ Geom.\ {\bf 1} (1992), 361--409.

\bibitem[Dem97]{Dem97} Demailly, J.-P., \textit{On the Ohsawa-Takegoshi-Manivel
$L^2$ extension theorem,} Proceedings of the Conference in honour of
the 85th birthday of Pierre Lelong, Paris, September 1997, \'ed.\
P.~Dolbeault, Progress in Mathematics, Birkh\"auser, Vol.~{\bf 188}
(2000), 47--82.

\bibitem[Dem15a]{Dem15a} 
Demailly, J.-P., 
\textit{On the cohomology of pseudoeffective line bundles,} 
Complex Geometry and Dynamics, The Abel Symposium 2013 held at Trondheim,
ed.\ John Erik Forn{\ae}ss, Marius Irgens, Erlend Forn{\ae}ss Wold,
Springer-Verlag 2015, 51--99.

\bibitem[Dem15b]{Dem15b}
Demailly, J.-P.,
\textit{Extension of holomorphic functions defined on non reduced analytic subvarieties, }
Advanced Lectures in Mathematics ALM35, 
the legacy of Bernhard Riemann after one hundred and fifty years,
Higher Education Press, Beijing-Boston, 2015, 191-222.

\bibitem[Dem-e-book]{Dem-e-book}
Demailly, J.-P.,
\textit {Complex analytic and differential geometry,} 
e-book available on the web page of the author,
{\tt https://www-fourier.ujf-grenoble.fr/\~{}demailly/manuscripts/agbook.pdf}

\bibitem[Dem-book]{Dem-book}
Demailly, J.-P.,  
\textit {Analytic methods in algebraic geometry,} 
Surveys of Modern Mathematics, {\bf{1}}. International Press, Somerville, 
Higher Education Press, Beijing, (2012).

\bibitem[DHP13]{DHP13} Demailly, J.-P., Hacon, Ch., P\u{a}un, M..
Extension theorems, Non-vanishing and the existence of good minimal models.
Acta Math.\ {\bf 210} (2013), 203--259.

\bibitem[DK01]{DK01} Demailly, J.-P., Koll\'ar, J.,
\textit{Semicontinuity of complex singularity exponents and K\"ahler-Einstein 
metrics on Fano orbifolds,} Ann.\ Ec.\ Norm.\ Sup.\ {\bf 34} (2001) 525--556.

\bibitem[DPS01]{DPS01}
Demailly, J.-P., Peternell, T., Schneider, M.,
\textit{Pseudo-effective line bundles on compact K\"ahler manifolds, }
International Journal of Math. {\bf{6}} (2001), 689--741. 

\bibitem[DF83]{DF83}
Donnelly, H., Fefferman, C.,
\textit{$L^2$-cohomology and index theorem for the Bergman metric,}
Ann.\ Math.\ {\bf 118} (1983) 593--618.

\bibitem[FK]{FK}
Folland, G.B., Kohn, J.J.,
\textit{The Neumann problem for the Cauchy-Riemann complex,}  
Annals of Mathematics Studies, No. {\bf{75}}. 
Princeton University Press, Princeton, 
N.J.; University of Tokyo Press, Tokyo, (1972). 

\bibitem[Fuj13]{Fuj13}
Fujino, O.,
\textit{A transcendental approach to Koll\'ar's injectivity theorem I\hspace{-.1em}I,}
J. Reine Angew. Math. {\bf{681}} (2013), 149--174. 

\bibitem[FM16]{FM16}
Fujino, O., Matsumura, S.-i.,
\textit{Injectivity theorem for pseudo-effective line bundles and its
applications,}
arXiv:1605.02284v1.

\bibitem[GZ15a]{GZ15a} Guan, Qi'an, Zhou, Xiangyu,
\textit{A solution of an $L^2$
extension problem with an optimal estimate and applications,}
Ann.\ of Math.\ (2) {\bf 181} (2015), no.~3, 1139--1208.

\bibitem[GZ15b]{GZ15b}
Guan Qi'an, Zhou, Xiangyu,
\textit{A proof of Demailly's strong openness conjecture,}
Annals of Math.\ {\bf 182} (2015), no.~2, 605--616.

\bibitem[Hie14]{Hie14}
Pham H.\ Hi$\rm{\hat{\underline{e}}}$p,  
\textit{The weighted log canonical threshold,}
C.~R.\ Math.\ Acad.\ Sci.\ Paris {\bf 352} (2014), 
no.~4, 283--288. 

\bibitem[H\"or65]{Hor65} H\"ormander, L., \textit{$L^2$ estimates and existence
theorems for the $\dbar$ operator,} Acta Math.\ {\bf 113}
(1965), 89--152.

\bibitem[H\"or66]{Hor66} H\"ormander, L., \textit{An introduction to Complex
Analysis in several variables,} 1st edition, Elsevier Science Pub.,
New York, 1966, 3rd revised edition, North-Holland Math.\ library, Vol
7, Amsterdam (1990).

\bibitem[Hos17a]{Hos17a} Hosono, Genki,
\textit{The optimal jet $L^2$ extension of Ohsawa-Takegoshi type,}
arXiv: math.CV/1706.08725

\bibitem[Hos17b]{Hos17b} Hosono, Genki,
\textit{On sharper estimates of Ohsawa-Takegoshi $L^2$-extension theorem,}
arXiv: math.CV/1708.08269.

\bibitem[Kod53a]{Kod53a} Kodaira, K., \textit{On cohomology groups of compact
analytic varieties with coefficients in some analytic faisceaux},
Proc.\ Nat.\ Acad.\ Sci.\ U.S.A.\ {\bf 39} (1953), 868--872.

\bibitem[Kod53b]{Kod53b} Kodaira, K., \textit{On a differential geometric
method in the theory of analytic stacks,} Proc.\ Nat.\ Acad.\ Sci.\
U.S.A.\ {\bf 39} (1953),
1268--1273.

\bibitem[Kod54]{Kod54} Kodaira, K., \textit{On K\"ahler varieties of
restricted type,} Ann.\ of Math.\ {\bf 60} (1954), 28--48.

\bibitem[Lem14]{Lem14}
Lempert, L.,
\textit{Modules of square integrable holomorphic germs,}
in: Analysis Meets Geometry: A Tribute to Mikael Passare,
Trends in Mathematics, Springer International Publishing, 2017, 311--333.

\bibitem[Man93]{Man93} Manivel, L., \textit{Un th\'eor\`eme de
prolongement $L^2$ de sections holomorphes d'un fibr\'e vectoriel,}
Math.\ Zeitschrift {\bf 212} (1993), 107--122.

\bibitem[Mat16]{Mat16}
Matsumura, S.-i.,
\textit{An injectivity theorem with multiplier ideal sheaves 
for higher direct images under K\"ahler morphisms,}
arXiv: math.CV/1607.05554v1. 

\bibitem[Mat17]{Mat17}
Matsumura, S.-i.,
\textit{An injectivity theorem with multiplier ideal sheaves 
of singular metrics with transcendental singularities,} arXiv:
math.CV/1308.2033v4,
J.~Algebraic Geom., DOI: \texttt{https://doi.org/10.1090/jag/687},
e-pub.\ August 17, 2017.

\bibitem[MV07]{MV07} McNeal, J., Varolin, D., \textit{Analytic inversion of adjunction: 
$L^2$ extension theorems with gain,} Ann.\ Inst.\ Fourier (Grenoble)
{\bf 57} (2007), no.~3, 703--718.

\bibitem[Nak55]{Nak55} Nakano, S., \textit{On complex analytic vector bundles},
J.~Math.\ Soc.\ Japan {\bf 7} (1955), 1--12.

\bibitem[Nad89]{Nad89} Nadel, A.M., \textit{Multiplier ideal sheaves and
K\"ahler-Einstein metrics of positive scalar curvature,} Proc.\ Nat.\
Acad.\ Sci.\ U.S.A. {\bf 86} (1989), 7299--7300~~~and~~
Annals of Math., {\bf 132} (1990), 549--596.

\bibitem[Nak73]{Nak73} Nakano, S., \textit{Vanishing theorems for weakly $1$-complete
manifolds,} Number Theory, Algebraic Geometry and Commutative
Algebra, in honor of Y.~Akizuki, Kinokuniya, Tokyo (1973), 169--179.

\bibitem[Nak74]{Nak74} Nakano, S., \textit{Vanishing theorems for weakly $1$-complete
manifolds~II,} Publ.\ R.I.M.S., Kyoto Univ.\ {\bf 10} (1974), 101--110.

\bibitem[Ohs88]{Ohs88} Ohsawa, T., \textit{On the extension of $L^2$ holomorphic
functions, II,} Publ. RIMS, Kyoto Univ.\ {\bf 24} (1988), 265--275.

\bibitem[Ohs94]{Ohs94} Ohsawa, T., \textit{On the extension of $L^2$ holomorphic
functions, IV$\,$: A new density concept,} Mabuchi, T.\ (ed.) et al.,
Geometry and analysis on complex manifolds. Festschrift for Professor
S.~Kobayashi's 60th birthday. Singapore: World Scientific,
(1994), 157--170.

\bibitem[Ohs95]{Ohs95} Ohsawa, T., \textit{On the extension of $L^2$ holomorphic
functions, III$\,$: negligible weights,} Math.\ Zeitschrift {\bf 219}
(1995), 215--225.

\bibitem[Ohs01]{Ohs01} Ohsawa, T., \textit{On the extension of $L^2$ holomorphic
functions, V$\,$: Effects of generalization,} Nagoya Math.\ J.\ {\bf 161}
(2001), 1--21, erratum: Nagoya Math.\ J.\ {\bf 163} (2001), 229.

\bibitem[Ohs03]{Ohs03} Ohsawa, T., \textit{On the extension of $L^2$ holomorphic
functions, VI$\,$: A limiting case,} Explorations in complex and Riemannian
geometry, Contemp.\ Math., 332, Amer.\ Math.\ Soc., Providence, RI,
2003, 235--239.

\bibitem[Ohs04]{Ohs04}
T. Ohsawa, 
\textit{On a curvature condition that implies a cohomology injectivity theorem 
of Koll\'ar-Skoda type,} 
Publ. Res. Inst. Math. Sci. {\bf{41}} (2005), no. 3, 565--577.

\bibitem[Ohs05]{Ohs05} Ohsawa, T., \textit{$L^2$ extension theorems -- backgrounds and
a new result,} Finite or infinite dimensional complex analysis and
applications, Kyushu University Press, Fukuoka, 2005, 261--274.

\bibitem[OT87]{OT87} Ohsawa, T., Takegoshi, K., \textit{On the extension of $L^2$
holomorphic functions,} Math.\ Zeitschrift {\bf 195} (1987), 197--204.

\bibitem[Pau07]{Pau07} P\u{a}un, M., 
\textit{Siu’s invariance of plurigenera: a one-tower proof,}
J.\ Differential Geom.\ {\bf 76} (2007) 485--493.

\bibitem[Pop05]{Pop05} Popovici, D., \textit{$L^2$ extension for jets of
holomorphic
sections of a Hermitian line bundle,} Nagoya Math.\ J.\ {\bf 180} (2005), 1--34.

\bibitem[Siu74]{Siu74} Siu, Y.-T.,
\textit{Analyticity of sets associated to Lelong
numbers and the extension of closed positive currents,}
Invent.\ Math.\ {\bf 27} (1974), 53--156.

\bibitem[Siu02]{Siu02}
Siu, Y.-T., \textit{Extension of twisted pluricanonical sections 
with plurisubharmonic weight and invariance of semipositively twisted 
plurigenera for manifolds not necessarily of general type}, Complex 
Geometry (G\"ottingen, 2000), Springer, Berlin, 2002, 223--277.

\bibitem[Var10]{Var10} Varolin, D., \textit{Three variations on a theme in
complex analytic geometry,} Analytic and Algebraic Geometry, IAS/Park City
Math.\ Ser.\ {\bf 17}, Amer.\ Math.\ Soc., (2010), 183--294.

\end{thebibliography}
\end{document}